\documentclass[reqno,12pt]{amsart}
\usepackage{amsmath, latexsym, amsfonts, amssymb, amsthm, amscd}

\numberwithin{equation}{section}

\newtheorem{theorem}{Theorem}[section]
\newtheorem{lemma}[theorem]{Lemma}
\newtheorem{proposition}[theorem]{Proposition}
\newtheorem{corollary}[theorem]{Corollary}
\newtheorem{rem}[theorem]{Remark}

\newtheorem{comment}[theorem]{Comment}

\renewcommand{\tilde}{\widetilde}          
\DeclareMathSymbol{\leqslant}{\mathalpha}{AMSa}{"36} 
\DeclareMathSymbol{\geqslant}{\mathalpha}{AMSa}{"3E} 
\DeclareMathSymbol{\eset}{\mathalpha}{AMSb}{"3F}     
\renewcommand{\leq}{\;\leqslant\;}                   
\renewcommand{\geq}{\;\geqslant\;}                   

\newcommand{\R}{\mathbb{R}}

\newcommand{\X}{\mathcal{X}}
\newcommand{\Y}{\mathcal{Y}}

\def\8{{\infty}}

\title[Hydrodynamic turbulence and intermittent random fields 
]{Hydrodynamic turbulence and intermittent random fields}
\author{}
\thanks{Partially
   supported by CNRS (UMR 7599
``Probabilit{\'e}s et Mod{\`e}les
Al{\'e}atoires'')}

\begin{document}

\maketitle
\begin{center}
{Raoul Robert \\
\footnotesize \noindent
 Institut Fourier, universit{\'e} Grenoble 1, UMR CNRS 5582, \\ 
 100, rue des Math{\'e}matiques, BP 74, 38402 Saint-Martin d'H{\`e}res cedex, France}

{\footnotesize \noindent e-mail: \texttt{Raoul.Robert@ujf-grenoble.fr}}

\bigskip

{Vincent Vargas \\
\footnotesize 
 CNRS, UMR 7534, F-75016 Paris, France \\
  Universit{\'e} Paris-Dauphine, Ceremade, F-75016 Paris, France} \\

 {\footnotesize \noindent e-mail: \texttt{vargas@ceremade.dauphine.fr}}
\end{center}

\begin{abstract}
In this article, we construct two families of multifractal random
vector fields with non symmetrical increments. We discuss the use of
such families to model the velocity field of turbulent flows.
\end{abstract}
\vspace{1cm}
\footnotesize
\noindent{\bf Short Title.} Hydrodynamic turbulence and intermittent random fields



\normalsize

\section{Introduction}
Roughly observed, some random phenomena seem scale
invariant. This is the case for the velocity field of turbulent flows
or the (logarithm of) evolution in time of the price of a
financial asset. However, a more precise empirical study of these
phenomena displays in fact a weakened form of scale invariance
commonly called multifractal scale invariance or intermittency (the
exponent which governs the power law scaling of the process or field
is no longer linear). An important question is therefore to construct
intermittent random fields which exhibit the observed characteristics.   

Following the work of Kolmogorov and Obukhov (\cite{cf:Kol}, \cite{cf:Obu})  on the energy
dissipation in turbulent flows, Mandelbrot introduced in \cite{cf:Man} a
"limit-lognormal" model to describe turbulent dissipation  or the
volatility of a financial asset. This model was rigorously defined and
studied in a mathematical framework by Kahane in \cite{cf:Kah}; more precisely,
Kahane constructed a random measure called Gaussian multiplicative
chaos. A natural extension of this work is to use Gaussian
multiplicative chaos to construct a field (or a process in the financial
case) which describes the whole phenomenon: the velocity field in
turbulent flows (the price of an asset on a financial market). This
extension was first performed by Mandelbrot himself who proposed to
model the price of a financial asset with a time changed Brownian
motion, the time change being random and independant of the Brownian
motion. In \cite{cf:BaDeMu}, the authors proposed for the time change to take the
primitive of multiplicative chaos: this gives the so called
multifractal random walk model (MRW) (Bacry and Muzy later generalized the
construction of the MRW model in \cite{cf:BaMuzy}). The obtained process accounts for
many observed properties of financial assets. 

The drawback of the above construction and of the MRW model 
is that the laws of the increments are symmetrical. In the case of
finance, this is in contradiction with the skewness property observed
for certain asset prices. In the case of turbulence, the laws of the
increments must be nonsymmetrical: it is a theoretical necessity and
stems from the dissipation of the kinetic energy (\cite{cf:Fri}). In light of these
observations, we are naturally led to construct
 random fields which generalize to any dimension such process and which
 present multifractal scale invariance as well as nonsymmetrical increments.          

We will answer a very natural question: how can one obtain a
family of multifractal fields with nonsymmetrical
increments by perturbing a given scale 
invariant Gaussian random field on $\R^d$?
Finally, in the last part we will mention the difficulties which arise
in trying to 
construct an incompressible multifractal velocity field that verifies
the $4/5$-law of Kolmogorov with positive dissipation.

\bigskip

\section{Notations and preliminary results}

\subsection{The underlying Gaussian field}

Let $dW_{0}(x)$ denote the Gaussian white noise on $\R^d$ and
$\varphi:\R^d \rightarrow [0,1]$ denote a $C^{\infty}$, radially
symmetric function worth $1$ for $|x| \leq 1$ and $0$ for $|x|>2$. 
We also introduce a fixed correlation scale $R >0$ and $\alpha$ a
number which satisfies  
\begin{equation}\label{eq:cond}
d/2< \alpha <d/2+1. 
\end{equation}
We define the Gaussian field ${\X}_{g}$ by
the following formula:

\begin{equation}\label{eq:gaussian}
{\X}_{g}(x)=\int_{\R^d}\varphi_{R}(x-y)\frac{x-y}{|x-y|^{d-\alpha+1}}dW_{0}(y),
\end{equation}
where we set the following notation:
\begin{equation*}
\varphi_{R}(x)=R^{d/2-\alpha}\varphi(\frac{x}{R}).
\end{equation*}
Using Kolmogorov's continuity criterion (see the standard book \cite{cf:ReYo}), it is easy to show that (\ref{eq:gaussian}) defines a homogeneous,
isotropic gaussian field which is almost surely H\"{o}lderian of order
$<\alpha-d/2$. Note that condition (\ref{eq:cond}) implies that the
integrand in (\ref{eq:gaussian}) is square integrable and the
$R^{d/2-\alpha}$ factor ensures that the field is dimensionless.

\subsubsection{Scaling property}
Let $e$ be a unitary vector and $\lambda>0$. We have the following
identity in law:
\begin{equation*}
{\X}_{g}(x+\lambda e)-{\X}_{g}(x)
\underset{(law)}{=}\int_{\R^d}(\frac{\varphi_{R}(y)y}{|y|^{d-\alpha+1}}
-\frac{\varphi_{R}(y-\lambda e)(y-\lambda e)}{|y-\lambda e|^{d-\alpha+1}})dW_{0}(y).
\end{equation*}
From the Gaussianity of the above law, we deduce that for all $q>0$,
there exists $c_{q}>0$ such that: 
\begin{equation*}
E(|{\X}_{g}(x+\lambda e)-{\X}_{g}(x)|^q)=\sigma_{\lambda e}^q c_{q},
\end{equation*}
with 
\begin{align*}
\sigma_{\lambda e}^2   & \; =   \int_{\R^d}(\frac{\varphi_{R}(y)y}{|y|^{d-\alpha+1}}
-\frac{\varphi_{R}(y-\lambda e)(y-\lambda e)}{|y-\lambda e|^{d-\alpha+1}})^2 dy  \\
& \; = \lambda^{2\alpha-d} \int_{\R^d}(\frac{\varphi_{R}(\lambda z)z}{|z|^{d-\alpha+1}}
-\frac{\varphi_{R}(\lambda z-\lambda e)(z-e)}{|z- e|^{d-\alpha+1}})^2 dz  \\
& \underset{ \lambda \rightarrow 0}{\sim}
(\frac{\lambda}{R})^{2\alpha-d}
\int_{\R^d}|\frac{z}{|z|^{d-\alpha+1}}-\frac{z-e}{|z-e|^{d-\alpha+1}}|^2 dz.
\end{align*}
We thus derive the following scaling 
\begin{equation*}
E(|{\X}_{g}(x+\lambda e)-{\X}_{g}(x)|^q)\underset{ \lambda \rightarrow 0}{\sim}(\frac{\lambda}{R})^{q(\alpha-d/2)}C_{q},
\end{equation*}
where the constant $C_{q}$ is independent of $e$. One says that
$({\X}_{g}(x))_{x \in \R^d}$ is at small scales monofractal with
scaling exponent $\alpha-d/2$.   

A homogeneous and isotropic field $({\X}(x))_{x \in \R^d}$ is multifractal if there exists a non
linear function $\zeta_{q}$ such that:
\begin{equation*}
E(|{\X}(x+\lambda e)-{\X}(x)|^q)\underset{ \lambda \rightarrow 0}{\sim}(\frac{\lambda}{R})^{\zeta_{q}}C_{q}.
\end{equation*}
We call $\zeta_{q}$ the structural function of the field $({\X}(x))_{x \in \R^d}$.

\subsection{Outline of the construction of multifractal vector fields from
  the field ${\X}_{g}$} 
Our construction is inspired by the work of Kahane in \cite{cf:Kah}. 
Let $\epsilon>0$ and $X^{\epsilon}(y)$ be a regular family of scalar Gaussian fields (not
necessarily independent of $dW_{0}$). We consider a family of fields
${\X}_{\epsilon}$ (with scalar components ${\X}_{\epsilon}^{j}$ in the canonical basis) defined by:
\begin{equation}\label{eq:cons}  
{\X}_{\epsilon}(x)=\int_{\R^d}\varphi_{R}(x-y)\frac{x-y}{|x-y|_{\epsilon}^{d-\alpha+1}}e^{X^{\epsilon}(y)-C_{\epsilon}}dW_{0}(y)
\end{equation}
($|x-y|_{\epsilon}$ is defined in the next subsection and is given by
a standard convolution). For an appropriate family $X^{\epsilon}$, we show that it is possible
to find constants $C_{\epsilon}$ such that ${\X}_{\epsilon}$ tends to
a non trivial field ${\X}$ (with scalar components ${\X}^{j}$ in the canonical basis) as $\epsilon$ tends to $0$.
If one chooses $X^{\epsilon}$ independent of $dW_{0}$, we will see
that this leads to a field ${\X}$ that extends the model introduced by
Bacry in \cite{cf:BaDeMu} and that has symmetrical increments. Thus, to obtain
nonsymmetrical increments, we must introduce correlation between
$X^{\epsilon}$ and $dW_{0}$. 

\subsection{Notations and construction of the family $X^{\epsilon}$}
Let $k^{R}$ be the function
\begin{equation*}
  k^R(x) =
  \begin{cases}
    \, \frac{1}{|x|^{d/2}} &\text{ for } |x| \leq R, \\
    \, 0 &\text{ otherwise}.
  \end{cases}
\end{equation*}

Let $\theta(x)$ be a $C^{\infty}$, non negative and radialy symmetrical
function with compact support in $|x| \leq 1$ such that 
\begin{equation*}
 \int_{\R^d}\theta(x)dx=1.
\end{equation*}
We define
$\theta^{\epsilon}=\frac{1}{\epsilon^d}\theta(\frac{.}{\epsilon})$ and 
the corresponding convolutions:
\begin{equation*}
k_{\epsilon}^R=\theta^{\epsilon}*k^{R}, \quad \quad
|.|_{\epsilon}=\theta^{\epsilon}*|.|. 
\end{equation*}
Let $\gamma$ be a strictly
positive parameter and $dW$ be a gaussian white noise on $\R^d$. We
consider the following gaussian field:
\begin{equation*} 
X^{\epsilon}(y)=\gamma \int_{\R^d}k_{\epsilon}^R(y-\sigma)dW(\sigma).
\end{equation*}
Its correlation kernel is given by:
\begin{equation*}
E(X^{\epsilon}(x) X^{\epsilon}(y))=\gamma^2\rho_{\epsilon/R}(\frac{x-y}{R}),
\end{equation*}
where $\rho=k^{1}*k^{1}$ and
$\rho_{\epsilon}=\theta^{\epsilon}*\theta^{\epsilon}*\rho$. One can
prove the following expansion
\begin{equation*}
\rho(x)=\omega_{d}\ln^{+}\frac{1}{|x|}+\phi(x),
\end{equation*}
where $\omega_{d}$ denotes the surface of the unit sphere in $\R^d$
and $\phi$ is a continuous function that vanishes for $|x| \geq 2$.
We will note $|.|_{*}=\text{inf}(1,|.|)$ and, with this definition, the
previous expansion is equivalent to:
\begin{equation*}
e^{\rho(x)}=\frac{e^{\phi(x)}}{|x|_{*}^{\omega_{d}}}.
\end{equation*}
One can also prove the following expansions with respect to $\epsilon$ for $\epsilon < R$:
\begin{equation}\label{eq:expansion1}
k_{\epsilon}^R(0)= \frac{C_{0}}{\epsilon^{d/2}}
\end{equation} 
with $C_{0}=\int_{|u| \leq 1} \frac{\theta(u)}{|u|^{d/2}}du$ and there exists
a constant $C_{1}$ such that 
\begin{equation}\label{eq:expansion2}
\rho_{\epsilon/R}(0)=\omega_{d}\ln\frac{R}{\epsilon} + C_{1} + o(\epsilon).
\end{equation} 
In the sequel, we will consider the case 
\begin{equation*} 
\gamma dW=\gamma_{0}(\epsilon)dW_{0}+\gamma_{1}dW_{1},
\end{equation*}
where $dW_{1}$ is a white noise independant of $dW_{0}$ and
$\gamma_{0}(\epsilon)$ is a function of $\epsilon$ that will be
defined later. Note that the integral in
formula (\ref{eq:cons}) has a meaning since $dW_{0}$ can be viewed as
a random distribution.

\subsection{Preliminary technical results}
We remind the following integration by parts formula for gaussian
vectors (cf. lemma 1.2.1 in \cite{cf:Nua}):

\begin{lemma}
Let $(g,g_{1}, \ldots, g_{n})$ be a centered gaussian vector and
$G:\R^{n} \rightarrow \R$ a $C^1$ function such that its partial
derivatives have at most exponential growth. Then we have:
\begin{equation}\label{eq:ipp}  
E(g G(g_{1}, \ldots, g_{n}))=\sum_{i=1}^{n}E(g
g_{i})E(\frac{\partial G}{\partial x_{i}}(g_{1}, \ldots, g_{n})).
\end{equation}
\end{lemma}

From the above formula, one can easily deduce by induction the
following lemma which will be frequently used in the sequel:

\begin{lemma}\label{lem:ind}
Let $l \in N^{*}$ be some positive integer and $(g,g_{1}, \ldots,
g_{2l})$ a centered gaussian vector. Then:
\begin{equation*}
E(g_{1} \ldots g_{2l}e^{g})=(\sum_{k=0}^{l}S_{k,l})e^{\frac{1}{2}E(g^2)},
\end{equation*}
where
\begin{equation*}
S_{k,l}=\sum_{\lbrace i_{1}, \ldots, i_{2k} \rbrace \subset \lbrace 1,
  \ldots, 2l \rbrace} \sum E(gg_{i_{1}}) \ldots
  E(gg_{i_{2k}})E(g_{i_{2k+1}}g_{i_{2k+2}}) \ldots E(g_{i_{2l-1}}g_{i_{2l}}),
\end{equation*}

where the second sum is taken over all partitions of $\lbrace 1,
  \ldots, 2l \rbrace \backslash \lbrace i_{1}, \ldots, i_{2k} \rbrace$ in subsets of two elements $\lbrace
  i_{2p+1},i_{2p+2} \rbrace$.

Similarly, we get the following formula: 

\begin{equation*}
E(g_{1} \ldots g_{2l+1}e^{g})=(\sum_{k=0}^{l}\tilde{S}_{k,l})e^{\frac{1}{2}E(g^2)},
\end{equation*}
where
\begin{equation*}
\tilde{S}_{k,l}=\sum_{\lbrace i_{1}, \ldots, i_{2k+1} \rbrace \subset \lbrace 1,
  \ldots, 2l+1 \rbrace} \sum E(gg_{i_{1}}) \ldots
  E(gg_{i_{2k+1}})E(g_{i_{2k+2}}g_{i_{2k+3}}) \ldots E(g_{i_{2l}}g_{i_{2l+1}}),
\end{equation*}
\end{lemma}

\begin{rem}
In  $S_{k,l}$ ($\tilde{S}_{k,l}$), the summation is made of
$\frac{2l!}{2k!2^{l-k}(l-k)!}$
($\frac{(2l+1)!}{(2k+1)!2^{l-k}(l-k)!}$) terms, number we will denote 
by $\alpha_{k,l}$ ($\tilde{\alpha}_{k,l}$). 
\end{rem}

We will also use the following lemma essentially due to Kahane (\cite{cf:Kah}).
\begin{lemma}\label{lem:Kah}
Let $(T,d)$ be a metric space and $\sigma$ a finite positive measure on $T$
equiped with the borelian $\sigma$-field induced by $d$. 

Let $q:T \times T: \rightarrow \R_{+}$ a symmetric application and $m$
a positive integer. Then we have the following inequalities:  
\begin{equation}\label{eq:kahp}
\int_{T^{2m}}e^{\sum_{1 \leq j < k \leq 2m}q(t_{j},t_{k})}d\sigma(t_{1})
\ldots d\sigma(t_{2m}) \leq \sigma(T) (\sup_{s \in T}\int_{T}e^{mq(t,s)}d\sigma(t))^{2m-1}, \end{equation}

\begin{align}
&  \int_{T^{2m+1}}e^{\sum_{1 \leq j < k \leq 2m+1}q(t_{j},t_{k})}d\sigma(t_{1})
\ldots d\sigma(t_{2m+1}) \nonumber \\  
&  \leq \sigma(T) \sup_{s,\tilde{s}}(\int_{T}
e^{q(\tilde{s},t)}d\sigma(t))(\int_{T}e^{mq(s,t)}e^{q(\tilde{s},t)}d\sigma(t))^{2m-1}.
\label{eq:kahi}
\end{align}
\end{lemma}
\proof
The proof of (\ref{eq:kahp}) can be found in \cite{cf:Kah}. 
Thus we just prove  how to derive inequality (\ref{eq:kahi}) from
 (\ref{eq:kahp}). By integrating with respect to the first $2m$
 variables and applying (\ref{eq:kahp}) with the measure $e^{q(t,t_{2m+1})}d\sigma(t)$, we get: 
 \begin{align*}
 &  \int_{T^{2m+1}}e^{\sum_{1 \leq j < k \leq 2m+1}q(t_{j},t_{k})}d\sigma(t_{1})
 \ldots d\sigma(t_{2m+1}) \\ 
 &  = \int_{T}d\sigma(t_{2m+1})\int_{T^{2m}}e^{\sum_{1 \leq j < k \leq
     2m}q(t_{j},t_{k})}\prod_{j=1}^{2m}e^{q(t_{j},t_{2m+1})}d\sigma(t_{1}) \ldots d\sigma(t_{2m}) \\ 
 &  \underset{(\ref{eq:kahp})}{\leq}\int_{T}d\sigma(t_{2m+1})(\int_{T}
 e^{q(t,t_{2m+1})}d\sigma(t))(\sup_{s}\int_{T}e^{mq(s,t)}e^{q(t,t_{2m+1})}d\sigma(t))^{2m-1} \\
 &  \leq \sigma(T) \sup_{s,\tilde{s}}(\int_{T}
 e^{q(\tilde{s},t)}d\sigma(t))(\int_{T}e^{mq(s,t)}e^{q(\tilde{s},t)}d\sigma(t))^{2m-1}.
 \end{align*}

\section{Construction of a four parameter family of multifractal,
  homogeneous, isotropic vector fields with non symmetrical increments }
\textbf{In this section, we will suppose that $d/2 < \alpha <
  (d/2+1) \wedge d$ and $\omega_{d}\gamma_{1}^{2}<d$}.
We consider the field $\X^{\epsilon}$ defined by formula
(\ref{eq:cons}) with 
\begin{equation*}
X^{\epsilon}(y)=\gamma_{0}(\epsilon)X_{0}^{\epsilon}(y)+\gamma_{1}X_{1}^{\epsilon}(y),
\end{equation*}
where 
\begin{equation*}
X_{i}^{\epsilon}(y)=\int_{\R^d}k_{\epsilon}^R(y-\sigma)dW_{i}(\sigma),
\qquad i=0,1.
\end{equation*}
We set also 
\begin{equation*}
C_{\epsilon}=((\gamma_{0}(\epsilon))^{2}+\gamma_{1}^2)\rho_{\epsilon/R}(0).
\end{equation*}
and
\begin{equation*}
\gamma_{0}(\epsilon)=\gamma_{0}^{*}(\frac{\epsilon}{R})^{\frac{d-\omega_{d}\gamma_{1}^2}{2}}.
\end{equation*}

Therefore, we introduce a slight correlation between $X^{\epsilon}$
and $dW_{0}$ ($\gamma_{0}(\epsilon)$ tends to $0$ as $\epsilon$ goes
to $0$).

\subsection{Multiplicative chaos in dimension $d$}
Multiplicative chaos or the "limit-lognormal" model introduced by
Mandelbrot is a generalization of the exponential of a gaussian
process. As mentioned in the introduction, it was defined rigorously
by Kahane in \cite{cf:Kah}. The construction of Kahane was based on the
theory of martingales and thus the generalized correlation kernel
(here $\rho(t-s)$) had to verify a
condition hard to verify practically (the $\sigma$-positivity
condition). Our construction is based on $L^2$-theory and can be
carried out without this condition. 

We will construct the multiplicative chaos associated to the
generalized correlation kernel $\rho(\frac{x-y}{R})$ defined in 2.3
and to some (positive) intermittency parameter $\gamma_{1}$ such that $\gamma_{1}^{2}\omega_{d}<d$.

Let $\epsilon$ be a positive number. Let $\mathcal{B}(\R^d)$ denote
the standard borelian $\sigma$-field; we want to consider the limit as $\epsilon$ goes to
$0$ of the random measures $Q^{\epsilon,\gamma_{1}}$ defined by:
\begin{align}
Q^{\epsilon,\gamma_{1}}(dy) & = e^{\gamma_{1}X_{1}^{\epsilon}(y)-\frac{\gamma_{1}^{2}}{2}E((X_{1}^{\epsilon}(y))^2)}dy
\nonumber \\
& =e^{\gamma_{1}X_{1}^{\epsilon}(y)-\frac{1}{2}\gamma_{1}^{2}\rho_{\epsilon/R}(0)}dy.\label{eq:chaoepsilon}
\end{align}

This leads us to state the following proposition:
\begin{proposition}[Multiplicative chaos of order $\gamma_{1}$]\label{prop:chaos}
There exists a positive random measure $Q^{\gamma_{1}}(dy)$
independent of the regularizing function $\theta$ such that:
\begin{enumerate}
\item
for all $A$ bounded in $\mathcal{B}(\R^d)$, $E(Q^{\gamma_{1}}(A))=|A|$.
\item
$Q^{\gamma_{1}}$ has almost surely no atoms.
\item
Almost surely, $Q^{\gamma_{1}}$ is singular with respect to the Lebesgue measure on all set $A$ (with positive measure). 
\end{enumerate}

If $q$ is some
positive integer and $f:\R^d \rightarrow \R$ a deterministic function
that satisfies the following condition:   
\begin{equation}\label{eq:condchao}
\int_{(\R^d)^{2q}}|f(y_{1})| \ldots |f(y_{2q})|
\prod_{1 \leq i < j\leq
 2q}\frac{1}{|\frac{y_{i}-y_{j}}{R}|_{*}^{\gamma_{1}^{2}\omega_{d}}}dy_{1}\ldots
   dy_{2q} < \infty,
\end{equation}
then we have the following convergence:
\begin{equation*}
\int_{\R^d}f(y)Q^{\epsilon,\gamma_{1}}(dy)\underset{\epsilon \to 0}{\overset{L^{2q}}{\rightarrow}}\int_{\R^d}f(y)Q^{\gamma_{1}}(dy).
\end{equation*}
We also have the following expression for the moments of $\int_{\R^d}f(y)Q^{\gamma_{1}}(dy)$:
\begin{equation}\label{eq:momchaos}
\forall k \leq 2q,  \quad
   E \left( (\int_{\R^d}f(y)Q^{\gamma_{1}}(dy))^{k} \right) =\int_{(\R^d)^k}f(y_{1})
   \ldots f(y_{k})
    \prod_{1 \leq i < j\leq
   k}\frac{e^{\gamma_{1}^{2}\phi(\frac{y_{i}-y_{j}}{R})}}{|\frac{y_{i}-y_{j}}{R}|_{*}^{\gamma_{1}^{2}\omega_{d}}}dy_{1}\ldots
   dy_{k}.
\end{equation}
Moreover, the above formula (\ref{eq:momchaos}) extends straightforwardly to the case of two functions 
$f,g$ and two intermittency parameters $\gamma_{1}, \gamma_{2}$ giving:
\begin{align*}
&  E \left( (\int_{\R^d}f(y)Q^{\gamma_{1}}(dy))^{k}(\int_{\R^d}g(y)Q^{\gamma_{2}}(dy))^{l} \right)\\
&=  \int_{(\R^d)^{k+l}}f(y_{1})
   \ldots f(y_{k})g(y_{k+1})
   \ldots g(y_{k+l})
    \prod_{1 \leq i < j\leq
   k}\frac{e^{\gamma_{1}^{2}\phi(\frac{y_{i}-y_{j}}{R})}}{|\frac{y_{i}-y_{j}}{R}|_{*}^{\gamma_{1}^{2}\omega_{d}}}  \times   \\
 &  \prod_{1 \leq i \leq
   k, \; j > k}\frac{e^{\gamma_{1}\gamma_{2}\phi(\frac{y_{i}-y_{j}}{R})}}{|\frac{y_{i}-y_{j}}{R}|_{*}^{\gamma_{1}\gamma_{2}\omega_{d}}} \prod_{k+1 \leq i < j\leq
   k+l}\frac{e^{\gamma_{2}^{2}\phi(\frac{y_{i}-y_{j}}{R})}}{|\frac{y_{i}-y_{j}}{R}|_{*}^{\gamma_{2}^{2}\omega_{d}}}  dy_{1}\ldots
   dy_{k+l}.
       \\
\end{align*}
We will call $Q^{\gamma_{1}}(dy)$ multiplicative chaos of order $\gamma_{1}$.
\end{proposition}

\proof
We first start by considering a positive integer $q$ and a function $f$ that satisfies the
corresponding integrability condition (\ref{eq:condchao}). 
Let $\epsilon, \epsilon'$ be two positive numbers. 
By using Fubini, we get for all $j \leq 2q$:
 \begin{align*}
  & E((\int_{\R^d}f(y)Q^{\epsilon,\gamma_{1}}(dy))^j(\int_{\R^d}f(y)Q^{\epsilon',\gamma_{1}}(dy))^{2q-j}) \\
    & = 
    e^{-\frac{j}{2}\gamma_{1}^{2}\rho_{\epsilon/R}(0)-\frac{2q-j}{2}\gamma_{1}^{2}\rho_{\epsilon'/R}(0)}\int_{(\R^d)^{2q}} f(y_{1})
   \ldots f(y_{2q})\times
   \\ 
 & \quad
    e^{\frac{1}{2}\gamma_{1}^{2}E((\sum_{i=1}^{j}X_{1}^{\epsilon}(y_{i})+\sum_{i=j+1}^{2q}X_{1}^{\epsilon'}(y_{i}))^2)}dy_{1}\ldots dy_{2q} \\
 & \underset{\epsilon,\epsilon' \to 0}{\rightarrow}\int_{(\R^d)^{2q}} f(y_{1})
   \ldots f(y_{2q}) \prod_{1 \leq
   i < j \leq 2q}\frac{e^{\gamma_{1}^{2}\phi(\frac{y_{i}-y_{j}}{R})}}{|\frac{y_{i}-y_{j}}{R}|_{*}^{\gamma_{1}^{2}\omega_{d}}}dy_{1}\ldots
   dy_{2q},
 \end{align*}
since $\rho_{\epsilon/R}(x) \rightarrow \rho(x)$ as $\epsilon$ goes to $0$.

From this, we deduce that:
 \begin{equation*}  
 E((\int_{\R^d}f(y)Q^{\epsilon,\gamma_{1}}(dy)-\int_{\R^d}f(y)Q^{\epsilon',\gamma_{1}}(dy))^{2q})\underset{\epsilon,\epsilon'
   \to 0}{\rightarrow}0
 \end{equation*}
 and therefore that $\int_{\R^d}f(y)Q^{\epsilon,\gamma_{1}}(dy)$ is a Cauchy sequence
 in $L^{2q}$ that converges to some random variable
 $\tilde{Q}^{\gamma_{1}}(f)$. For $k \leq 2q$, the moment
 $E((\tilde{Q}^{\gamma_{1}}(f))^k)$ is the limit as $\epsilon$ goes to
 $0$ of $E((Q^{\epsilon,\gamma_{1}}(f))^k)$; from this one can deduce
 that the moments of $\tilde{Q}^{\gamma_{1}}(f)$ are given by formula (\ref{eq:momchaos}).

For any bounded set $A$ in $\mathcal{B}(\R^d)$, consider $f=1_{A}$ and $q=1$. Since $\gamma_{1}^{2}\omega_{d}<d$, we deduce from lemma \ref{lem:Kah} that the integrability condition (\ref{eq:condchao}) is satisfied. Thus it follows from the proof above that  
$Q^{\epsilon,\gamma_{1}}(A)$ converges in $L^2$ to some random
variable $\tilde{Q}^{\gamma_{1}}(A)$. This defines a family of random
variables (indexed by the bounded Borelian sets) that satisfies the
following properties:
\begin{enumerate}
\item
For all disjoint and bounded sets $A_{1}$, $A_{2}$ in $\mathcal{B}(\R^d)$, 
\begin{equation*}
\tilde{Q}^{\gamma_{1}}(A_{1} \cup A_{2})=\tilde{Q}^{\gamma_{1}}(A_{1})+\tilde{Q}^{\gamma_{1}}(A_{2})
\quad \quad a.s.
\end{equation*}
\item
For any bounded sequence $(A_{n})_{n \geq 1}$ decreasing to $\varnothing$:  
\begin{equation*}
\tilde{Q}^{\gamma_{1}}(A_{n})\underset{n \to \infty}{\rightarrow}0
\quad a.s.
\end{equation*}
\end{enumerate}

By theorem 6.1.VI. in \cite{cf:DaVe}, there exists a random measure
$Q^{\gamma_{1}}$ such that for all bounded $A$ in $\mathcal{B}(\R^d)$
we have: 
\begin{equation*} 
Q^{\gamma_{1}}(A)=\tilde{Q}^{\gamma_{1}}(A) \quad a.s.
\end{equation*}

Finally, one can easily show that the limit random variable
$\tilde{Q}^{\gamma_{1}}(f)$ is almost surely equal to $\int_{\R^d}f(y)Q^{\gamma_{1}}(dy)$.
\qed

\subsection{Convergence of $\X_{\epsilon}$ towards a field $\X$.}

\vspace{1em}

\textbf{}

 \textbf{In the sequel, $(e_{j})_{j}$ will denote the
canonical basis (whereas $(e^{j})_{j}$ denotes the components of a
vector $e$).} In this subsection, we will prove the following proposition: 

\begin{proposition}\label{prop:conv}
Let $\alpha$ be such that $d/2 < \alpha < (d/2+1) \wedge d$ and $\gamma_{1}$ such
that $2\gamma_{1}^{2}\omega_{d}<\alpha-d/2$. There exists a field
$(\X(x))_{x \in \R^d}$ such that for all $k$ and $x_{1}, \ldots, x_{k}
\in \R^d$ the following convergence in law holds: 
\begin{equation}\label{eq:conv} 
 (\X_{\epsilon}(x_{1}), \ldots, \X_{\epsilon}(x_{k}))
 \underset{\epsilon \to 0}{\Rightarrow} (\X(x_{1}), \ldots, \X(x_{k})).
\end{equation}

Let $l$ be an
integer such that one of the following conditions hold:
\begin{enumerate}
\item
$l$ is even and $l\gamma_{1}^{2}\omega_{d}<\alpha-d/2$.
\item
$l$ is odd and $(l+1)\gamma_{1}^{2}\omega_{d}<\alpha-d/2$.
\end{enumerate}
Let $F_{R}^{j}$ be defined by $F_{R}^{j}(y)=\varphi_{R}(y)\frac{y^{j}}{|y|^{d-\alpha+1}}$.
Then there exists $C$ such that, for all $x$ in $\R^d$, the random variables $\X^{j}(x)$ have a moment
of order $2l$ given by the following expression:
\begin{eqnarray}
& & E((\X^{j}(x))^{2l}) = \sum_{k=0}^l\alpha_{k,l}C^{2k}\int_{(\R^d)^{k+l}}
F_{R}^{j}(y_{1}) \ldots F_{R}^{j}(y_{2k}) (F_{R}^{j}(y_{2k+1}))^{2}
\ldots (F_{R}^{j}(y_{k+l}))^{2}
  \nonumber   \\ 
& & \prod_{1 \leq i < j \leq
  2k}\frac{e^{\gamma_{1}^{2}\phi(\frac{y_{i}-y_{j}}{R})}}{|\frac{y_{i}-y_{j}}{R}|_{*}^{\gamma_{1}^{2}\omega_{d}}}\prod_{\underset{
  j > 2k}{1 \leq i \leq
  2k}}\frac{e^{2\gamma_{1}^{2}\phi(\frac{y_{i}-y_{j}}{R})}}{|\frac{y_{i}-y_{j}}{R}|_{*}^{2\gamma_{1}^{2}\omega_{d}}}
  \prod_{2k+1 \leq i < j \leq
  k+l}\frac{e^{4\gamma_{1}^{2}\phi(\frac{y_{i}-y_{j}}{R})}}{|\frac{y_{i}-y_{j}}{R}|_{*}^{4\gamma_{1}^{2}\omega_{d}}}dy_{1}
  \ldots dy_{k+l}. \nonumber \\
& & \label{eq:mom1}
\end{eqnarray}

We also have: 
\begin{eqnarray}
& & E((\X^{j}(x+h)-\X^{j}(x))^{2l})=\sum_{k=0}^l\alpha_{k,l}C^{2k}\int_{(\R^d)^{k+l}}
(F_{R}^{j}(y_{1})-F_{R}^{j}(y_{1}-h))
  \ldots \nonumber \\
& & (F_{R}^{j}(y_{2k})-F_{R}^{j}(y_{2k}-h))
(F_{R}^{j}(y_{2k+1})-F_{R}^{j}(y_{2k+1}-h))^2 \ldots
  (F_{R}^{j}(y_{k+l})-F_{R}^{j}(y_{k+l}-h))^2 \nonumber  \\ 
& & \prod_{1 \leq i < j \leq
  2k}\frac{e^{\gamma_{1}^{2}\phi(\frac{y_{i}-y_{j}}{R})}}{|\frac{y_{i}-y_{j}}{R}|_{*}^{\gamma_{1}^{2}\omega_{d}}}\prod_{\underset{
  j > 2k}{1 \leq i \leq
  2k}}\frac{e^{2\gamma_{1}^{2}\phi(\frac{y_{i}-y_{j}}{R})}}{|\frac{y_{i}-y_{j}}{R}|_{*}^{2\gamma_{1}^{2}\omega_{d}}}
  \prod_{2k+1 \leq i < j \leq
  k+l}\frac{e^{4\gamma_{1}^{2}\phi(\frac{y_{i}-y_{j}}{R})}}{|\frac{y_{i}-y_{j}}{R}|_{*}^{4\gamma_{1}^{2}\omega_{d}}}dy_{1}
  \ldots dy_{k+l}. \label{eq:mom2} 
\end{eqnarray}

\end{proposition}

\proof

Let $\gamma_{1}$ be such that $2\gamma_{1}^{2}\omega_{d}<\alpha-d/2$. We set 
\begin{equation*}  
C=\frac{\gamma_{0}^{*}C_{0}e^{-1/2\gamma_{1}^{2}C_{1}}}{R^{d/2}}.
\end{equation*}
and define two auxiliary fields $\Y_{\epsilon}, \mathcal{Z}_{\epsilon}$
by the following expressions:
\begin{equation}\label{eq:champaux1}
\Y_{\epsilon}(x)=\int_{\R^d}\varphi_{R}(x-y)\frac{x-y}{|x-y|^{d-\alpha+1}}Q^{\epsilon,\gamma_{1}}(dy)
\end{equation}
and
\begin{equation}\label{eq:champaux2}
\mathcal{Z}_{\epsilon}(x)=\int_{\R^d}\varphi_{R}(x-y)\frac{x-y}{|x-y|^{d-\alpha+1}}
e^{\gamma_{1}X_{1}^{\epsilon}(y)-C_{\epsilon}}dW_{0}(y).
\end{equation}
Note that $\mathcal{Z}_{\epsilon}(x)$ exists since $X_{1}^{\epsilon}$
and $dW_{0}$ are independent with:
\begin{equation}
E(\int_{\R^d}\frac{\varphi_{R}(x-y)^2}{|x-y|^{2(d-\alpha)}}
e^{2\gamma_{1}X_{1}^{\epsilon}(y)-2C_{\epsilon}}dy) < \infty.
\end{equation}

We can compute, for all $x$ in $\R^d$,
$E(|\X_{\epsilon}(x)-C\Y_{\epsilon}(x)-\mathcal{Z}_{\epsilon}(x)|^2)$
(cf. the more complicated computations in the proof of proposition \ref{prop:tight}) 
and derive the following limit:
\begin{equation*}
\X_{\epsilon}(x)-(C\Y_{\epsilon}(x)+\mathcal{Z}_{\epsilon}(x))\underset{\epsilon
\to 0}{\overset{L^2}{\rightarrow}}0.
\end{equation*}
Thus, we must show that the finite dimensional distributions of the
field $C\Y_{\epsilon}+\mathcal{Z}_{\epsilon}$ converge in law. 
Let $k$ be some positive integer and $x_{1}, \ldots, x_{k}$ points in
$\R^d$. For all $\xi=(\xi_{1}, \ldots, \xi_{k})$ in $(\R^d)^k$, we
compute the characteristic function of
$(C\Y_{\epsilon}(x_{1})+\mathcal{Z}_{\epsilon}(x_{1}),\ldots,C\Y_{\epsilon}(x_{k})+\mathcal{Z}_{\epsilon}(x_{k}))$:
\begin{equation*}
\mathcal{C}_{\epsilon}(\xi)=E(e^{i\sum_{j=1}^{k}\xi_{j}.(C\Y_{\epsilon}(x_{j})+\mathcal{Z}_{\epsilon}(x_{j}))}).
\end{equation*}

By conditioning on the field generated by the white noise $dW_{1}$, we get:

\begin{align*}
\mathcal{C}_{\epsilon}(\xi) & =E(e^{iC\sum_{j=1}^{k}\xi_{j}.\Y_{\epsilon}(x_{j})}e^{-\frac{1}{2}\int(\sum_{j=1}^{k}\xi_{j}.F_{R}(x_{j}-y))^{2}e^{2\gamma_{1}X_{1}^{\epsilon}(y)-2C_{\epsilon}}dy}) \\
&=E(e^{iC\sum_{j=1}^{k}\int\xi_{j}.F_{R}(x_{j}-y)Q^{\epsilon,\gamma_{1}}(dy)}e^{-\frac{1}{2}e^{-2\gamma_{0}^{2}\rho_{\epsilon/R}(0)}\int(\sum_{j=1}^{k}\xi_{j}.F_{R}(x_{j}-y))^{2}Q^{\epsilon,2\gamma_{1}}(dy)})
\end{align*}

Now, using proposition (\ref{prop:chaos}), we have:   
\begin{equation*}
\sum_{j=1}^{k}\int\xi_{j}.F_{R}(x_{j}-y)Q^{\epsilon,\gamma_{1}}(dy)\underset{L^{2}}{\rightarrow}\sum_{j=1}^{k}\int\xi_{j}.F_{R}(x_{j}-y)Q^{\gamma_{1}}(dy),
\end{equation*}

\begin{equation*}
\int(\sum_{j=1}^{k}\xi_{j}.F_{R}(x_{j}-y))^{2}Q^{\epsilon,2\gamma_{1}}(dy)\underset{L^{2}}{\rightarrow}\int(\sum_{j=1}^{k}\xi_{j}.F_{R}(x_{j}-y))^{2}Q^{2\gamma_{1}}(dy),
\end{equation*}

from where:
\begin{equation*}
\mathcal{C}_{\epsilon}(\xi)\underset{\epsilon \to 0}{\rightarrow}\mathcal{C}(\xi)=E(e^{iC\sum_{j=1}^{k}\int\xi_{j}.F_{R}(x_{j}-y)Q^{\gamma_{1}}(dy)}e^{-\frac{1}{2}\int(\sum_{j=1}^{k}\xi_{j}.F_{R}(x_{j}-y))^{2}Q^{2\gamma_{1}}(dy)}).
\end{equation*}

Thus, by applying Levy's theorem, we conclude that the finite
dimensional distributions of the field
$C\Y_{\epsilon}+\mathcal{Z}_{\epsilon}$ converge in law to those of a
field $\X$ whose finite dimensional distributions are given by: 
\begin{equation*}
E(e^{i\sum_{j=1}^{k}\xi_{j}.\X(x_{j})})=\mathcal{C}(\xi).
\end{equation*}
Suppose that $l$ is a positive integer that satisfies the condition of
the proposition. For all $\xi$ in $\R^{d}$, we have:
  
\begin{equation*}
E(e^{i\xi.\X(x)})=E(e^{-iC\int\xi.F_{R}(y)Q^{\gamma_{1}}(dy)}e^{-\frac{1}{2}\int(\xi.F_{R}(y))^{2}Q^{2\gamma_{1}}(dy)}).
\end{equation*}
We derive expression (\ref{eq:mom1}) by computing
$\frac{\partial^{2l}}{(\partial\xi^{j})^{2l}}E(e^{i\xi.\X(x)})|_{\xi=0}$
using proposition  \ref{prop:chaos}. We derive (\ref{eq:mom2}) similarly.

\subsection{Scaling of $\X$}
The purpose of this subsection is 
to show that the field
$(\X(x))_{x \in \R^d}$ satisfies the multifractal scaling relation
(this is what propositions \ref{prop:scap} and \ref{prop:scai}
below assert).

We first state two preliminary lemmas we will use in the rest of the paper.

\begin{lemma}
Let $\delta$ be some real number such that $0 \leq \delta <
\alpha$ and $\delta \not = \alpha-1$. There exists $C=C(\delta)$ such that we have the following
inequality for $|h| \leq R$:
\begin{equation}\label{eq:estuni1}
\sup_{x \in
  \R^d}\int_{\R^d}|F_{R}(y)-F_{R}(y-h)|\frac{1}{|\frac{x-y}{R}|_{*}^{\delta}}dy
  \leq R^{d/2} C |\frac{h}{R}|^{(\alpha -\delta)\wedge{1}}. 
\end{equation}
\end{lemma}

\proof
By homogenity, we suppose that $R=1$ and for simplicity, we suppose
$d\geq 2$. Since $\frac{1}{|x|_{*}^{\delta}} \leq  \frac{1}{|x|^{\delta}}+1$ and the right hand side of 
(\ref{eq:estuni1}) increases with $\delta$, we have to show that
for $\delta \in [0,\alpha[$ and $|h| \leq 1$:
\begin{equation*}
\sup_{x \in
  \R^d}\int_{\R^d}|F_{1}(y)-F_{1}(y-h)|
\frac{1}{|x-y|^{\delta}}dy \leq C|h|^{(\alpha -\delta)\wedge{1}}. 
\end{equation*}
Indeed, this would imply that for $|h| \leq 1$:
\begin{align*}
\sup_{x \in
  \R^d}\int_{\R^d}|F_{1}(y)-F_{1}(y-h)|
\frac{1}{|x-y|_{*}^{\delta}}dy   &  \leq C|h|^{(\alpha -\delta)\wedge{1}}+ C|h|^{\alpha \wedge{1}} \\
& \leq 2 C |h|^{(\alpha -\delta)\wedge{1}}.  \\  
\end{align*}

There exists $C$ such that for all $y$ and $h$, we have:
\begin{equation}\label{eq:lip}
|\varphi(y-h)-\varphi(y)| \leq C|h|\text{  and  }
 \varphi(y) \leq C1_{|y| \leq
  2}.
\end{equation}

We set 
\begin{equation*}
I(x)=\int_{\R^d}|F_{1}(y)-F_{1}(y-h)|
\frac{1}{|x-y|^{\delta}}dy.
\end{equation*}

Therefore we get
\begin{align}
I(x) & \leq  C|h|\int_{|y| \leq 3}\frac{1}{|y|^{d-\alpha}}\frac{1}{|x-y|^{\delta}}dy \\ 
&  \quad +C\int_{|y|\leq 3}|\frac{y}{|y|^{d-\alpha+1}}-\frac{y-h}{|y-h|^{d-\alpha+1}}|\frac{1}{|x-y|^{\delta}}dy
  \nonumber \\
&  \leq  C|h|+C\int_{|y| \leq
  3}|\frac{y}{|y|^{d-\alpha+1}}-\frac{y-h}{|y-h|^{d-\alpha+1}}|\frac{1}{|x-y|^{\delta}}dy,
\label{eq:pasdn}
\end{align}
where we denote by $C$ different constants.

\noindent{\textit{First case}: $\delta<\alpha-1$.}

Plugging inequality 
\begin{equation*}
|\frac{y}{|y|^{d-\alpha+1}}-\frac{y-h}{|y-h|^{d-\alpha+1}}| \leq \frac{(d-\alpha+1)|h|}{|y-h|^{d-\alpha+1}\wedge|y|^{d-\alpha+1}}  
\end{equation*}
in (\ref{eq:pasdn}), we get 
\begin{equation*}
I(x) \leq C|h|\int_{|y| \leq 3}\frac{1}{|y-h|^{d-\alpha+1}\wedge|y|^{d-\alpha+1}}\frac{1}{|x-y|^{\delta}}dy.
\end{equation*}
We have:
\begin{align*}
& \int_{|y| \leq 3}\frac{1}{|y-h|^{d-\alpha+1}\wedge|y|^{d-\alpha+1}}\frac{1}{|x-y|^{\delta}}dy \\
& \leq  \int_{|y| \leq 3}\frac{1}{|y-h|^{d-\alpha+1}}\frac{1}{|x-y|^{\delta}}dy +
 \int_{|y| \leq 3}\frac{1}{|y|^{d-\alpha+1}}\frac{1}{|x-y|^{\delta}}dy \\
 & \leq 2 \sup_{x}\int_{|y| \leq 4}\frac{1}{|y|^{d-\alpha+1}}\frac{1}{|x-y|^{\delta}}dy,
\end{align*}

which concludes the proof.

\noindent{\textit{Second case}: $\delta>\alpha-1$.}

By the change of variable
$y=|h|u$ and setting $h=|h|e$ with $|e|=1$, we get: 
\begin{align*}
&  \int_{|y| \leq
  3}|\frac{y-h}{|y-h|^{d-\alpha+1}}-\frac{y}{|y|^{d-\alpha+1}}|\frac{1}{|x-y|^{\delta}}dy \\
&  = |h|^{\alpha-\delta}\int_{|u| \leq
  \frac{3}{|h|}}|\frac{u-e}{|u-e|^{d-\alpha+1}}-\frac{u}{|u|^{d-\alpha+1}}|\frac{1}{|x/|h|-u|^{\delta}}du \\
&  \leq |h|^{\alpha-\delta} \sup_{a \in
  \R^d}\int_{\R^d}|\frac{u-e}{|u-e|^{d-\alpha+1}}-\frac{u}{|u|^{d-\alpha+1}}|\frac{1}{|a-u|^{\delta}}du.
\end{align*}

\qed

\begin{lemma}
Let $\delta$ be some real number such that $0 \leq \delta <
2\alpha-d$. There exists $C=C(\delta)$ such that we have the following
inequality for $|h| \leq R$:
\begin{equation}\label{eq:estuni2}
\sup_{x \in
  \R^d}\int_{\R^d}|F_{R}(y)-F_{R}(y-h)|^{2}\frac{1}{|\frac{x-y}{R}|_{*}^{\delta}}dy
  \leq C|\frac{h}{R}|^{2\alpha-d-\delta}. 
\end{equation}
\end{lemma}

\proof

As in the proof above, we can replace $|.|_{*}$ by $|.|$ and suppose
that $R=1$; thus we
have to show inequality (\ref{eq:estuni2}) with $J(x)$ where we set:
\begin{equation*}
J(x)=\int_{\R^d}|F_{R}(y)-F_{R}(y-h)|^{2}
\frac{1}{|x-y|^{\delta}}dy.
\end{equation*}
Using inequality (\ref{eq:lip}), we get 
\begin{align}
  J(x) & \leq C|h|^{2}\int_{|y| \leq
  3}\frac{1}{|y|^{2(d-\alpha)}}\frac{1}{|x-y|^{\delta}}dy
  \nonumber \\ 
& \quad + C\int_{|y|
  \leq
  3}|\frac{y-h}{|y-h|^{d-\alpha+1}}-\frac{y}{|y|^{d-\alpha+1}}|^{2}\frac{1}{|x-y|^{\delta}}dy
  \nonumber \\
 &   \leq C|h|^{2}+C\int_{|y|
  \leq
  3}|\frac{y-h}{|y-h|^{d-\alpha+1}}-\frac{y}{|y|^{d-\alpha+1}}|^{2}\frac{1}{|x-y|^{\delta}}dy. \label{eq:pasdn1} 
\end{align}
Since $2>2\alpha-d-\delta$, we only have to consider the second
term in inequality (\ref{eq:pasdn1}). By the change of variable
$y=|h|u$ and setting $h=|h|e$ with $|e|=1$, we get:
\begin{align*}
&  \int_{|y| \leq
  3}|\frac{y-h}{|y-h|^{d-\alpha+1}}-\frac{y}{|y|^{d-\alpha+1}}|^{2}\frac{1}{|x-y|^{\delta}}dy \\
&  = |h|^{2\alpha-d-\delta}\int_{|u| \leq
  \frac{3}{|h|}}|\frac{u-e}{|u-e|^{d-\alpha+1}}-\frac{u}{|u|^{d-\alpha+1}}|^{2}\frac{1}{|x/|h|-u|^{\delta}}du \\
&  \leq |h|^{2\alpha-d-\delta} \sup_{a \in \R^d}\int_{\R^d}|\frac{u-e}{|u-e|^{d-\alpha+1}}-\frac{u}{|u|^{d-\alpha+1}}|^{2}\frac{1}{|a-u|^{\delta}}du.\end{align*}

\qed

\begin{proposition}\label{prop:scap}{(Scaling along the even integers)}

Let $l$ be an
integer such that one of the following conditions hold:

\begin{enumerate}
\item
$l$ is even and $l\gamma_{1}^{2}\omega_{d}<\alpha-d/2$
\item
$l$ is odd and $(l+1)\gamma_{1}^{2}\omega_{d}<\alpha-d/2$.
\end{enumerate}
Let $e$ be a unit vector ($|e|=1$). Then there exists $C_{l}^{j}(e) >
0$ such that the following scaling relation holds:
\begin{equation}\label{eq:sca}
E((\X^{j}(x+\lambda e)-\X^{j}(x))^{2l})\underset{ \lambda \to 0}{\sim}C_{l}^{j}(e)(\frac{\lambda}{R})^{\zeta_{2l}},
\end{equation}
where we have
\begin{equation}\label{eq:zeta}
\zeta_{2l}=l(2\alpha-d)-2\gamma_{1}^2\omega_{d}l(l-1).
\end{equation}
\end{proposition}

\proof
For simplicity, we will suppose that $l$ is even and that
$l\gamma_{1}^{2}\omega_{d}<\alpha-d/2$. We introduce the following notation:  
\begin{equation*}
f_{h}(y)=F_{R}^{j}(y)-F_{R}^{j}(y-h).
\end{equation*}
We shall see that the scaling at small scale of the sum
(\ref{eq:mom2}) is given by the term $k=0$. Indeed for all $k \geq 1$
let us consider the integral 

\begin{eqnarray}
& & \int_{(\R^d)^{k+l}}
f_{h}(y_{1})
  \ldots f_{h}(y_{2k})
(f_{h}(y_{2k+1}))^2 \ldots
  (f_{h}(y_{k+l}))^2 \nonumber  
 \prod_{1 \leq i < j \leq
  2k}\frac{e^{\gamma_{1}^{2}\phi(\frac{y_{i}-y_{j}}{R})}}{|\frac{y_{i}-y_{j}}{R}|_{*}^{\gamma_{1}^{2}\omega_{d}}} \times \\
& & \prod_{\underset{
  j > 2k}{1 \leq i \leq
  2k}}\frac{e^{2\gamma_{1}^{2}\phi(\frac{y_{i}-y_{j}}{R})}}{|\frac{y_{i}-y_{j}}{R}|_{*}^{2\gamma_{1}^{2}\omega_{d}}}
  \prod_{2k+1 \leq i < j \leq
  k+l}\frac{e^{4\gamma_{1}^{2}\phi(\frac{y_{i}-y_{j}}{R})}}{|\frac{y_{i}-y_{j}}{R}|_{*}^{4\gamma_{1}^{2}\omega_{d}}}dy_{1}
  \ldots dy_{k+l}. \nonumber \\
& & \leq I_{k}J_{k,l},
\end{eqnarray}
where we set 
\begin{align*}
 I_{k}= & \sup_{y_{2k+1}, \ldots,
  y_{k+l}}\int_{(\R^d)^{2k}} |f_{h}(y_{1})| \ldots |f_{h}(y_{2k})|\prod_{\underset{
  j > 2k}{1 \leq i \leq
  2k}}\frac{e^{2\gamma_{1}^{2}\phi(\frac{y_{i}-y_{j}}{R})}}{|\frac{y_{i}-y_{j}}{R}|_{*}^{2\gamma_{1}^{2}\omega_{d}}}
  \times \\  
& \prod_{1 \leq i < j \leq
  2k}\frac{e^{\gamma_{1}^{2}\phi(\frac{y_{i}-y_{j}}{R})}}{|\frac{y_{i}-y_{j}}{R}|_{*}^{\gamma_{1}^{2}\omega_{d}}}dy_{1}
  \ldots dy_{2k}.
\end{align*}
and 
\begin{equation*}
 J_{k,l}= \int_{(\R^d)^{l-k}}(f_{h}(y_{2k+1}))^2 \ldots (f_{h}(y_{k+l}))^2
 \prod_{2k+1 \leq i < j \leq
 k+l}\frac{e^{4\gamma_{1}^{2}\phi(\frac{y_{i}-y_{j}}{R})}}{|\frac{y_{i}-y_{j}}{R}|_{*}^{4\gamma_{1}^{2}\omega_{d}}}dy_{2k+1}
 \ldots dy_{k+l}.
\end{equation*}
By using the estimates (\ref{eq:estuni1}),(\ref{eq:estuni2}) and the
inequalities (\ref{eq:kahp}), (\ref{eq:kahi}), one can show that for
all $k \geq 1$, we have   
\begin{equation*}
I_{k}J_{k,l} \leq C R^{dk}|\frac{h}{R}|^{c_{k,l}},
\end{equation*}
with 
\begin{equation*}
c_{k,l}=(\alpha-2(l-k)\gamma_{1}^{2}\omega_{d})\wedge 1 +((\alpha-(2l-k)\gamma_{1}^{2}\omega_{d})\wedge 1)(2k-1)+(2\alpha-d)(l-k)-2\gamma_{1}^{2}\omega_{d}(l-k)(l-k-1). 
\end{equation*}

If $\alpha-2(l-k)\gamma_{1}^{2}\omega_{d}<1$, then $c_{k,l}=\zeta_{2l}+k(d-\gamma_{1}^{2}\omega_{d})$;
If $\alpha-2(l-k)\gamma_{1}^{2}\omega_{d} \geq 1$ and
$\alpha-(2l-k)\gamma_{1}^{2}\omega_{d} < 1$, then
$c_{k,l}=\zeta_{2l}+1-\alpha+dk+(2l-3k)\gamma_{1}^{2}\omega_{d}$; 
otherwise $c_{k,l}=2k+(2\alpha-d)(l-k)-2\gamma_{1}^{2}\omega_{d}(l-k)(l-k-1)$.
In all cases, it is easy to show that $c_{k,l}>\zeta_{2l}$ under the
conditions of the proposition.

Finally, we study the term where $k=0$. We get for $h=\lambda e$ with $|e|=1$: 
\begin{align*}
&  \int_{(\R^d)^{l}}(f_{h}(y_{1}))^2 \ldots
  (f_{h}(y_{l}))^2\prod_{1 \leq i < n \leq
  l}\frac{e^{4\gamma_{1}^{2}\phi(\frac{y_{i}-y_{n}}{R})}}{|\frac{y_{i}-y_{n}}{R}|_{*}^{4\gamma_{1}^{2}\omega_{d}}}dy_{1}
  \ldots dy_{l}.  \\
&  \underset{y_{i}=\lambda
  u_{i}}{=}(\frac{\lambda}{R})^{l(2\alpha-d)}\int_{(\R^d)^{l}}(\varphi(\frac{\lambda}{R}(u_{1}-e))\frac{u_{1}^{j}-e^{j}}{|u_{1}-e|^{d-\alpha+1}}-\varphi(\frac{\lambda}{R}u_{1})\frac{u_{1}^{j}}{|u_{1}|^{d-\alpha+1}})^2 \ldots \times \\
 & \quad (\varphi(\frac{\lambda}{R}(u_{l}-e))\frac{u_{l}^{j}-e^{j}}{|u_{l}-e|^{d-\alpha+1}}-\varphi(\frac{\lambda}{R}u_{l})\frac{u_{l}^{j}}{|u_{l}|^{d-\alpha+1}})^2  
  \prod_{1 \leq i < n \leq
  l}\frac{e^{4\gamma_{1}^{2}\phi(\frac{\lambda(u_{i}-u_{n})}{R})}}{|\frac{\lambda(u_{i}-u_{n})}{R}|_{*}^{4\gamma_{1}^{2}\omega_{d}}}du_{1}
  \ldots du_{l}.     \\
&  \underset{\lambda \to
  0}{\sim}e^{2l(l-1)\gamma_{1}^{2}\phi(0)}(\frac{\lambda}{R})^{\zeta_{2l}}\int_{(\R^d)^{l}}(\frac{u_{1}^{j}-e^{j}}{|u_{1}-e|^{d-\alpha+1}}-\frac{u_{1}^{j}}{|u_{1}|^{d-\alpha+1}})^2 \ldots
  (\frac{u_{l}^{j}-e^{j}}{|u_{l}-e|^{d-\alpha+1}}-\frac{u_{l}^{j}}{|u_{l}|^{d-\alpha+1}})^2
  \times \\ 
& \quad  \prod_{1 \leq i < n \leq
  l}\frac{1}{|u_{i}-u_{n}|^{4\gamma_{1}^{2}\omega_{d}}}du_{1}
  \ldots du_{l}, 
\end{align*}
and inequality (\ref{eq:kahp}) shows that this integral is finite when $l\gamma_{1}^{2}\omega_{d}<\alpha-d/2$.

\qed

In the next proposition, we state the scaling relations of $\X$ along
the odd integers. We define $I_{l}^{j}(e)$ by:
\begin{align*}
I_{l}^{j}(e) & =\int_{(\R^d)^{l}}(\frac{u_{1}^{j}-e^{j}}{|u_{1}-e|^{d-\alpha+1}}-\frac{u_{1}^{j}}{|u_{1}|^{d-\alpha+1}})^2 \ldots
  (\frac{u_{l}^{j}-e^{j}}{|u_{l}-e|^{d-\alpha+1}}-\frac{u_{l}^{j}}{|u_{l}|^{d-\alpha+1}})^2
  \times \\ 
& \quad  \prod_{1 \leq i < j \leq
  l}\frac{1}{|u_{i}-u_{j}|^{4\gamma_{1}^{2}\omega_{d}}}du_{1}
  \ldots du_{l}. 
\end{align*}

\begin{proposition}\label{prop:scai}{(Scaling along the odd integers)}
Let $l$ be an
integer satisfying the conditions in proposition
\ref{prop:scap} and $1+2l\gamma_{1}^{2}\omega_{d}<\alpha$. 

Let $e$ be a unit vector ($|e|=1$).
Then we have the following scaling relation:

\begin{equation}\label{eq:scai}
E((\X^{j}(x+\lambda e)-\X^{j}(x))^{2l+1})\underset{ \lambda \to 0}{\sim}\Sigma_{l}^{j}(e)(\frac{\lambda}{R})^{\tilde{\zeta}_{2l+1}},
\end{equation}
where we have
\begin{equation}\label{eq:zetai}
\tilde{\zeta}_{2l+1}=l(2\alpha-d)-2\gamma_{1}^2\omega_{d}l(l-1)+1,
\end{equation}
and
\begin{equation*}
\Sigma_{l}^{j}(e)=\gamma_{0}^{*}C(l,\gamma_{1})I_{l}^{j}(e)e^{j},
\quad C(l,\gamma_{1}) >0.
\end{equation*}
\end{proposition}

\proof

As in proposition \ref{prop:conv}, setting
$C=\frac{\gamma_{0}^{*}C_{0}e^{-1/2\gamma_{1}^{2}C_{1}}}{R^{d/2}}$, it
is possible to show that:
 \begin{align*}
&  E((\X^{j}(x+h)-\X^{j}(x))^{2l+1}) \\
& =\sum_{k=0}^l\tilde{\alpha}_{k,l}C^{2k+1}\int_{(\R^d)^{k+l+1}}
f_{h}(y_{1})
  \ldots f_{h}(y_{2k+1})
(f_{h}(y_{2k+2}))^2 \ldots
  (f_{h}(y_{k+l+1}))^2 \times  \\ 
& \prod_{1 \leq i < j \leq
  2k+1}\frac{e^{\gamma_{1}^{2}\phi(\frac{y_{i}-y_{j}}{R})}}{|\frac{y_{i}-y_{j}}{R}|_{*}^{\gamma_{1}^{2}\omega_{d}}}\prod_{\underset{
  j > 2k+1}{1 \leq i \leq
  2k+1}}\frac{e^{2\gamma_{1}^{2}\phi(\frac{y_{i}-y_{j}}{R})}}{|\frac{y_{i}-y_{j}}{R}|_{*}^{2\gamma_{1}^{2}\omega_{d}}}
  \prod_{2k+2 \leq i < j \leq
  k+l+1}\frac{e^{4\gamma_{1}^{2}\phi(\frac{y_{i}-y_{j}}{R})}}{|\frac{y_{i}-y_{j}}{R}|_{*}^{4\gamma_{1}^{2}\omega_{d}}}dy_{1}
  \ldots dy_{k+l+1}, 
\end{align*}
where, as usual, we set:  
\begin{equation*}
f_{h}(y)=F_{R}^{j}(y)-F_{R}^{j}(y-h).
\end{equation*}
Similarly to proposition \ref{prop:scap}, to get the main contribution
as $|h|$ goes to $0$, we examine the term $k=0$.
We introduce $\mathcal{I}$: 
\begin{align*}
\mathcal{I} &= \int_{(\R^d)^{l+1}}f_{h}(y_{1})(f_{h}(y_{2}))^2 \ldots
   (f_{h}(y_{l+1}))^2\prod_{j \geq
   2}e^{2\gamma_{1}^{2}\rho(\frac{y_{1}-y_{j}}{R})}
   \times \\
 & \prod_{2 \leq i < j \leq
   l+1}\frac{e^{4\gamma_{1}^{2}\phi(\frac{y_{i}-y_{j}}{R})}}{|\frac{y_{i}-y_{j}}{R}|_{*}^{4\gamma_{1}^{2}\omega_{d}}}dy_{1}
   \ldots dy_{l+1}.
\end{align*}
Putting $h=\lambda e$ ($|e|=1$), $y_{1}=Ru_{1}$, $y_{i}=\lambda u_{i}$ ($i
\geq 2$), we get:
\begin{equation*}
\mathcal{I} \underset{\lambda \to 0} {\sim}
R^{d/2}e^{2l(l-1)\gamma_{1}^{2}\phi(0)}(\frac{\lambda}{R})^{\tilde{\zeta}_{2l+1}}I_{l}^{j}(e)
e.\int \nabla \psi^{j}(y)e^{2l\gamma_{1}^{2}\rho(y)}dy,
\end{equation*}
where $\psi^{j}(y)=\varphi(y)\frac{y^{j}}{|y|^{d-\alpha+1}}$.
Now a calculation gives:
\begin{equation*}
\int \nabla \psi^{j}(y)e^{2l\gamma_{1}^{2}\rho(y)}dy=-\frac{2l\gamma_{1}^{2}\omega_{d}}{d}(\int_{0}^{\infty}r^{\alpha-1}\varphi(r)e^{2l\gamma_{1}^{2}\rho(r)}\frac{d\rho}{dr}dr)e_{j},
\end{equation*}
the last integral being negative since $\rho(r)$ is a strictly
decreasing function on the interval $]0,2[$: the result follows. 
Notice that the condition $1+2l\gamma_{1}^{2}<\alpha$ implies that this
integral is finite while the other conditions imply that $I_{l}^{j}(e)$ is finite.
 
\qed

From proposition \ref{prop:scai}, we deduce readily that for $\lambda$
small the law of $\X(x+\lambda e)-\X(x)$ is nonsymmetrical (for
$\gamma_{0}^{*} \not = 0$). Indeed, by isotropy, we have :

\begin{equation*}
(\X(x+\lambda e)-\X(x)).e\underset{law}{=}\X^{j}(x+\lambda e_{j})-\X^{j}(x)
\end{equation*}
and
$\Sigma_{l}^{j}(e_{j})=\gamma_{0}^{*}C(l,\gamma_{1})I_{l}^{j}(e_{j}) >
0$.

\subsection{Tightness of $\X_{\epsilon}$ and regularity of $\X$.}

In this section, we prove that the convergence in law of $\X_{\epsilon}$ towards $\X$ given by proposition \ref{prop:conv} holds in a functional sense and that the field $\X$ is locally H\"{o}lderian. The straight way to do so is to prove the tightness of the sequence $\X_{\epsilon}$ by means of a Kolmogorov estimate (cf. chapter 13 of \cite{cf:ReYo}).

\begin{proposition}\label{prop:tight}{(Tightness)}
Let $l$ be some positive integer that satisfies the condition of
proposition \ref{prop:scap} and $\gamma$ a positive parameter such
that $\gamma_{1}^{2} < \gamma^{2}$. Then there exists $\epsilon_{0}>0$
and $C$ independent of
$\epsilon$ such that for $\epsilon < \epsilon_{0}$ and $|h| \leq R$:
 
\begin{equation}\label{eq:tight}
 \forall x, \qquad E((\X_{\epsilon}(x+h)-\X_{\epsilon}(x))^{2l}) \leq C|h|^{l(2\alpha-d)-2\gamma^2\omega_{d}l(l-1)},
\end{equation}
and  
\begin{equation}\label{eq:tight'}
E((\X_{\epsilon}(0))^{2l})\leq C.
\end{equation}

\end{proposition}

\proof

We only prove (\ref{eq:tight}) (the proof of \ref{eq:tight'} is similar). 
We are going to compute the moment
 \begin{equation*}
 E((\X_{\epsilon}^{j}(x+h)-\X_{\epsilon}^{j}(x))^{2l})=
E((\int_{\R^d}f_{\epsilon,h}(y)e^{X^{\epsilon}(y)-C_{\epsilon}}dW_{0}(y))^{2l})
 \end{equation*}
 where we set: 
 \begin{equation*}
 f_{\epsilon,h}(y)=\frac{\varphi_{R}(y)y^{j}}{|y|_{\epsilon}^{d-\alpha+1}}-\frac{\varphi_{R}(y-h)(y^{j}-h^{j})}{|y-h|_{\epsilon}^{d-\alpha+1}}.
 \end{equation*}
 We get: 
 \begin{equation}\label{eq:mom3}
 E((\int_{\R^d}f_{\epsilon,h}(y)e^{X^{\epsilon}(y)-C_{\epsilon}}dW_{0}(y))^{2l})
 =e^{-2l C_{\epsilon}}\int f_{\epsilon,h}(y_{1}) \ldots
 f_{\epsilon,h}(y_{2l})E(e^{\hat{X}^{\epsilon}}dW_{0}(y_{1}) \ldots
 dW_{0}(y_{2l})),
 \end{equation}
 where
 \begin{equation*} 
 \hat{X}^{\epsilon}=X^{\epsilon}(y_{1})+ \ldots + X^{\epsilon}(y_{2l}).
 \end{equation*}
 The rest of the computation can be performed rigorously by
 regularizing the white noise $dW_{0}$, using lemma \ref{lem:ind} and going to
 the limit. It is easy to see that we obtain the same result by
 introducing the following formal rules:
 \begin{equation}\label{eq:rule1}
 E(dW_{0}(y)dW_{0}(y'))=\delta_{y-y'}dy
 \end{equation}
 and
 \begin{equation}\label{eq:rule2}
 E(dW_{0}(y)X^{\epsilon}(y'))=\gamma_{0}(\epsilon)k_{\epsilon}^R(y'-y)dy
 \end{equation}
As a consequence of lemma \ref{lem:ind}, $E(e^{\hat{X}^{\epsilon}}dW_{0}(y_{1}) \ldots
dW_{0}(y_{2l}))$ is the sum of terms of the form
\begin{equation}\label{eq:term}
E(dW_{0}(y_{1})\hat{X}^{\epsilon}) \ldots
E(dW_{0}(y_{k})\hat{X}^{\epsilon})E(dW_{0}(y_{k+1})dW_{0}(y_{k+2}))
\ldots E(dW_{0}(y_{q-1})dW_{0}(y_{2l})) e^{\frac{1}{2}E((\hat{X}^{\epsilon})^{2})}.
\end{equation}
We will compute the limit of each one of these terms. By using
(\ref{eq:rule2}), we get 
\begin{align*}
E(dW_{0}(y_{k})\hat{X}^{\epsilon}) & = 
\gamma_{0}(\epsilon)(\sum_{i=1}^{2l} k_{\epsilon}^R(y_{i}-y_{k}))dy_{k} \\
& =  \gamma_{0}(\epsilon) k_{\epsilon}^R(0)(1+Q_{k}^{\epsilon})dy_{k} \\
\end{align*}
where 
\begin{equation*}
Q_{k}^{\epsilon}=\frac{1}{k_{\epsilon}^R(0)}(\sum_{i \not = k}k_{\epsilon}^R(y_{i}-y_{k})). 
\end{equation*}
We also have from the definition of $X^{\epsilon}$:
\begin{equation*}
e^{\frac{1}{2}E((\hat{X}^{\epsilon})^{2})}=e^{(l\rho_{\epsilon/R}(0)+\sum_{i
      < j}\rho_{\epsilon/R}(\frac{y_{i}-y_{j}}{R}))(\gamma_{0}(\epsilon)^2+\gamma_{1}^2)}.
\end{equation*}

 By using lemma \ref{lem:ind}, expression (\ref{eq:mom3}) and the rules
 above, we get:
 \begin{eqnarray*}
 & E((\int_{\R^d}f_{\epsilon,h}(y)e^{X^{\epsilon}(y)-C_{\epsilon}}dW_{0}(y))^{2l})
  = 
 \sum_{k=0}^{l}\alpha_{k,l}(\gamma_{0}(\epsilon))^{2k}(k_{\epsilon}^R(0))^{2k}e^{(2l-k)((\gamma_{0}(\epsilon))^{2}+\gamma_{1}^2)\rho_{\epsilon/R}(0)-2l
   C_{\epsilon}} & \\ 
  & \int_{(\R^{d})^{k+l}} f_{\epsilon,h}(y_{1}) \ldots
 f_{\epsilon,h}(y_{2k})(f_{\epsilon,h}(y_{2k+1}))^2 \ldots
 (f_{\epsilon,h}(y_{k+l}))^2
 \prod_{i=1}^{2k}(1+Q_{i,k,l}^{\epsilon})e^{S_{k,l}^{\epsilon}}dy_{1}
 \ldots dy_{k+l} &  \\ 
 \end{eqnarray*}
 where 
 \begin{equation*}
 Q_{i,k,l}^{\epsilon}=\frac{1}{k_{\epsilon}^R(0)}(\sum_{\underset{j
     \not = i}{1 \leq j \leq 2k}}k_{\epsilon}^R(y_{i}-y_{j})+2\sum_{j>2k}k_{\epsilon}^R(y_{i}-y_{j}))
 \end{equation*}
 and 
 \begin{eqnarray*}
  S_{k,l}^{\epsilon}& = & ((\gamma_{0}(\epsilon))^{2}+\gamma_{1}^2)(\sum_{1
   \leq i < j \leq 2k}\rho_{\epsilon/R}(\frac{y_{i}-y_{j}}{R})+2\sum_{1
   \leq i \leq 2k} \sum_{j >
   2k}\rho_{\epsilon/R}(\frac{y_{i}-y_{j}}{R})  \\ 
 & & +4\sum_{2k+1 \leq i < j \leq k+l}\rho_{\epsilon/R}(\frac{y_{i}-y_{j}}{R})). 
 \end{eqnarray*}
 We first take care of the normalizing constant outside each
 integral:
 \begin{equation*}
 (\gamma_{0}(\epsilon)k_{\epsilon}^R(0)e^{-1/2((\gamma_{0}(\epsilon))^{2}+\gamma_{1}^2)\rho_{\epsilon/R}(0)})^{2k}e^{2l((\gamma_{0}(\epsilon))^{2}+\gamma_{1}^2)\rho_{\epsilon/R}(0)-2l
   C_{\epsilon}}. 
 \end{equation*}
 By the choice of $C_{\epsilon}$, we have $e^{2l((\gamma_{0}(\epsilon))^{2}+\gamma_{1}^2)\rho_{\epsilon/R}(0)-2l
   C_{\epsilon}}=1$. Using expansions  (\ref{eq:expansion1}) and
   (\ref{eq:expansion2}), we derive the following limit:
 \begin{equation*}
 \gamma_{0}(\epsilon)k_{\epsilon}^R(0)e^{-1/2((\gamma_{0}(\epsilon))^{2}+\gamma_{1}^2)\rho_{\epsilon/R}(0)}\underset{
   \epsilon \to 0}{\rightarrow} \frac{\gamma_{0}^{*}C_{0}e^{-1/2\gamma_{1}^{2}C_{1}}}{R^{d/2}}.
 \end{equation*}
 In conclusion, the constant outside the integral of term $k$ in the
 above sum converges to $\alpha_{k,l}(\frac{\gamma_{0}^{*}C_{0}e^{-1/2\gamma_{1}^{2}C_{1}}}{R^{d/2}})^{2k}$.

Let $\gamma$ be such that $\gamma_{1}^2<\gamma^{2}$. One can choose $\epsilon_{0}>0$ such that
$\gamma_{0}(\epsilon_{0}))^{2}+\gamma_{1}^2<\gamma^{2}$. Using the fact that, for all $y$, $\rho_{\epsilon/R}(y/R)
\leq \omega_{d}\ln^{+}\frac{R}{|y|}+C$ with $C$ independent of
$\epsilon$, we get:
 \begin{align}
  e^{S_{k,l}^{\epsilon}}  \leq C \nonumber
  & \prod_{1 \leq i < j \leq
   2k}\frac{1}{|\frac{y_{i}-y_{j}}{R}|_{*}^{\gamma^{2}\omega_{d}}}\prod_{\underset{ j > 2k}{1 \leq i \leq 2k}}\frac{1}{|\frac{y_{i}-y_{j}}{R}|_{*}^{2\gamma^{2}\omega_{d}}} \times \\ 
 &  \prod_{2k+1 \leq i < j \leq
   k+l}\frac{1}{|\frac{y_{i}-y_{j}}{R}|_{*}^{4\gamma^{2}\omega_{d}}}.\label{eq:leb}
 \end{align}
Finally, we conclude by using the fact that $|Q_{i,k,l}^{\epsilon}|$
is bounded by a constant independent of $\epsilon$, inequality (\ref{eq:estuni1}) and
(\ref{eq:estuni2}) similarly as in the proof of proposition \ref{prop:scap}. 

\qed
\begin{corollary}
One can easily deduce
from this that for $\gamma_{1}^{2}$ sufficiently small, by
Kolmogorov's compacity theorem (\cite{cf:ReYo}), $\X_{\epsilon}$ tends to $\X$ in the
functional sense and that $\X$ is locally H\"{o}lderian.
\end{corollary}

\begin{comment}
Starting with a two parameter ($R,\alpha$) monofractal Gaussian field, we constructed a four
parameter ($R,\alpha, \gamma_{1}, \gamma_{0}^{*}$) multifractal field with
nonsymmetrical increments. This family has it's
own interest. As we shall see in the next section, this family is too
restricted to take into account all the constraints needed for a
satisfactory model of turbulent flows. 


In the case where $\gamma_{0}^{*}=0$, we obtain symmetrical random
fields which extend to higher dimensions the model introduced in \cite{cf:BaDeMu}.
\end{comment}

In the next section, we will study a multifractal field  which is not
in this family but that can be seen as a limit case where
$\gamma_{1}=0$ and  $\gamma_{0}$ is constant (independent of
$\epsilon$). As we will see, this family will be compatible with the $4/5$-law.

\section{A step towards a model of the velocity field of turbulent flows}
An acceptable solution to the problem of hydrodynamical turbulence in
dimension $3$ would be to construct a random velocity field $U$
solution to the dynamics (Euler or Navier Stokes typically) that is stationnary,
incompressible, space-homogeneous, isotropic and that satisfies the
main statistical properties of the velocity field of turbulent
flows. Two main properties are:
\begin{enumerate}     
\item
The $4/5$-law of Kolmogorov that links the energy dissipation of the
turbulent flow to the statistics of the increments of the
velocity. This law is widely accepted since it is the only one that
can be proven with the dynamics (\cite{cf:DuRonon}, \cite{cf:Fri}, \cite{cf:Ro}). More
precisely, this law states:   
\begin{equation}\label{eq:4/5}
E((U(x+\xi)-U(x).\frac{\xi}{|\xi|})^{3})=-\frac{4}{5}D|\xi|.
\end{equation}
In the above formula, $D$ denotes the average dissipation of the
kinetic energy per unit mass in the fluid.
\begin{rem}
To obtain this law, it is sufficient to suppose that the field $U$ is space
homogeneous and isotropic. 
\end{rem}
\item
The intermittency of the field $U$:
\begin{equation}\label{eq:inter}
E((U(x+\xi)-U(x).\frac{\xi}{|\xi|})^{q})\underset{|\xi| \to 0}{\sim}
C_{q} |\xi|^{\zeta_{q}},
\end{equation}
where $\zeta_{q}$ is a well known concave structure function (cf. \cite{cf:Fri}).
\end{enumerate}

It is a very challenging task to construct a field with all the
aforementioned properties, especially because this field must be
invariant by the Euler or Navier-Stokes equation. 

Nevertheless, one can in the first place forget the invariance
by the dynamics and simply try to construct a field that satisfies the
other properties. The $4/5$-law shows that the nonsymmetry of the
increments is an essential feature.
Let us consider the family $\X$ constructed in the previous section ($d=3$).
By proposition \ref{prop:scai}, we have:
\begin{equation*}
E(((\X(x+\lambda e)-\X(x)).e)^{3})\underset{\lambda \to 0}{\sim}
C_{3}(\frac{\lambda}{R})^{\tilde{\zeta}_{3}}, \quad C_{3} >0,
\end{equation*}
with $\tilde{\zeta}_{3}=2\alpha-2$.
To satisfy the $4/5$ law one should have $\tilde{\zeta}_{3}=1$, which
gives $\alpha=3/2$. This is incompatible with the constraint $3/2<
\alpha <5/2$.
Thus we have now to modify the family $\X$ to reach the limit case
$\tilde{\zeta}_{3}=1$. In this aim, we will construct a new (three
parameter) family $\X_{0}$ corresponding to the limit case
$\gamma_{1}=0$, $\gamma_{0}$ constant $>0$.

\subsection{Construction of the field $\X_{0}$}

 In this section, we only outline the main steps of the construction of
$\X_{0}$.
The field $\X_{0,\epsilon}$ is given by formula (\ref{eq:cons})
where $X^{\epsilon}$ is now defined by:
\begin{equation*}
X^{\epsilon}(y)=\gamma_{0}\int_{\R^{d}}k_{\epsilon}^{R}(y-\sigma)dW_{0}(\sigma).
\end{equation*}
We suppose that $\alpha$ is in the interval $]0,1[$.
We choose the normalizing constant $C_{\epsilon}$ such that:
\begin{equation*}
\gamma_{0}k_{\epsilon}^{R}(0)e^{-C_{\epsilon}+\frac{1}{2}\gamma_{0}^{2}\rho_{\epsilon/R}(0)}=1.
\end{equation*}
We start by stating a lemma we will use in the proof of the
proposition below: 
\begin{lemma}
let $\delta$ be some real number different from $d$. Then there
exists $C=C(\delta)>0$ with: 
\begin{equation}\label{eq:epsilon}
\int_{|u| \leq R}\frac{du}{|u|_{\epsilon}^{\delta}} \leq C
 \epsilon^{(d-\delta) \wedge 0}.
\end{equation}
\end{lemma}

\proof
We suppose $\delta > d$, the other case being obvious. We have:
\begin{align*}
 \int_{|u| \leq R}\frac{du}{|u|_{\epsilon}^{\delta}} 
& \underset{u=\epsilon  \tilde{u}}  {=} \epsilon^{d-\delta} \int_{|u| \leq R/
  \epsilon}\frac{d\tilde{u}}{(\int_{|v| \leq
    1}\theta(v)|v+\tilde{u}|dv)^{\delta}} \\
& \leq \epsilon^{d-\delta} \int_{\R^{d}}\frac{d\tilde{u}}{(\int_{|v| \leq 1}\theta(v)|v+\tilde{u}|dv)^{\delta}}.
\end{align*}
\qed

We can now state the following proposition:

\begin{proposition}\label{prop:scalingimp}
Let $l$ be an integer $\geq 1$ and $\gamma_{0}$ such that:
\begin{enumerate}
 \item
 $\gamma_{0}^{2}\omega_{d}<\alpha$ if $l=1$.
 \item
$(2l-3/2)\gamma_{0}^{2}\omega_{d}<\alpha\wedge d/2$ if $l>1$.
 \end{enumerate}
Then for all $x$, $\X_{0,\epsilon}(x)$ converges in $L^{2l}$ to a random
vector $\X_{0}(x)$ (i.e. $E((\X_{0,\epsilon}(x)-\X_{0}(x))^{2l})
\rightarrow 0$). The random vector field $\X_{0}(x)$ satisfies the
following scaling:   
For $e$ ($|e|=1$), and for all $q \leq 2l$: 
\begin{equation}\label{eq:scaling}
E((\X_{0}^{j}(x+\lambda e)-\X_{0}^{j}(x))^{q})\underset{\lambda \to
  0}{\sim}C_{q}^{j}(e)(\frac{\lambda}{R})^{\zeta_{q}},
\end{equation}
where $\zeta_{q}=q\alpha-\frac{1}{2}q(q-1)\gamma_{0}^{2}\omega_{d}$ and
\begin{equation}\label{eq:truccc}
C_{q}^{j}(e)=R^{qd/2}e^{\frac{q(q-1)}{2}\gamma_{0}^{2}\phi(0)}\int_{(\R^{d})^{q}}\prod_{1 \leq i < j \leq
  q}\frac{1}{|u_{i}-u_{j}|^{\gamma_{0}^{2}\omega_{d}}}\prod_{1
   \leq i \leq q}(\frac{u_{i}^{j}}{|u_{i}|^{d-\alpha+1}}-\frac{u_{i}^{j}-e^{j}}{|u_{i}-e|^{d-\alpha+1}})du_{1} \ldots du_{q}.
\end{equation}

\end{proposition}

\proof
We will first prove that:
\begin{equation}\label{eq:moment}
 E((\X_{0,\epsilon}^{j}(x))^{2l})\underset{\epsilon \to 0}{\rightarrow}
 \int_{(\R^{d})^{2l}}\prod_{1 \leq i < j \leq
  2l}\frac{e^{\gamma_{0}^{2}\phi(\frac{y_{i}-y_{j}}{R})}}{|\frac{y_{i}-y_{j}}{R}|_{*}^{\gamma_{0}^{2}\omega_{d}}}\prod_{1
  \leq i \leq 2l}\frac{\varphi_{R}(y_{i})y_{i}^{j}}{|y_{i}|^{d-\alpha+1}}dy_{1} \ldots dy_{2l}.
\end{equation}
We remind that the right hand side of the above limit exists by lemma \ref{lem:Kah}.
In order to prove the above relation, we develop
$E((\X_{0,\epsilon}^{j}(x))^{2l})$ in $l+1$ terms similarly as in the
proof of proposition \ref{prop:tight}; then, using formula
(\ref{eq:expansion1}) and the fact that, for all $y$, $\rho_{\epsilon/R}(y)
\leq \omega_{d}\ln\frac{R}{\epsilon}+C$, we are led to show that, for
all $k \leq l-1$, we have the following convergence: 
\begin{align*}
&
\epsilon^{(l-k)(d-\gamma_{0}^{2}\omega_{d})}\epsilon^{-2(l-k)(l-k-1)\gamma_{0}^{2}\omega_{d}}\epsilon^{-4k(l-k)\gamma_{0}^{2}\omega_{d}}\int_{(\R^{d})^{k+l}}\frac{\varphi_{R}(y_{1})y_{1}^{j}}{|y_{1}|_{\epsilon}^{d-\alpha+1}}
\ldots \frac{\varphi_{R}(y_{2k})y_{2k}^{j}}{|y_{2k}|_{\epsilon}^{d-\alpha+1}}
\times \\
& \frac{(\varphi_{R}(y_{2k+1})y_{2k+1}^{j})^2}{|y_{2k+1}|_{\epsilon}^{2(d-\alpha+1)}} \ldots
  \frac{(\varphi_{R}(y_{k+l})y_{k+l}^{j})^2}{|y_{k+l}|_{\epsilon}^{2(d-\alpha+1)}} \nonumber   
 \prod_{1 \leq i < j \leq
  2k}\frac{e^{\gamma_{0}^{2}\phi(\frac{y_{i}-y_{j}}{R})}}{|\frac{y_{i}-y_{j}}{R}|_{*}^{\gamma_{0}^{2}\omega_{d}}}dy_{1} \ldots dy_{k+l}
\underset{\epsilon \to 0}{\rightarrow} 0.
\end{align*}
We apply inequality (\ref{eq:epsilon}) and
obtain (if $\alpha=d/2$, one can work with $\alpha-\eta$ for $\eta>0$
sufficiently small) : 
\begin{equation*}
\int_{\R^{d}}\frac{\varphi_{R}(y)^2}{|y|_{\epsilon}^{2(d-\alpha)}}dy
\leq C \epsilon^{(2\alpha-d)\wedge 0}
\end{equation*}
Therefore the above convergence to $0$ amounts to showing that, for
all $k \leq l-1$, we have the following inequality:
\begin{equation*}
d+(2\alpha-d)\wedge 0-\gamma_{0}^{2}\omega_{d}>2(l-k-1)\gamma_{0}^{2}\omega_{d}+4k\gamma_{0}^{2}\omega_{d}.
\end{equation*}
This is equivalent to $(2l-\frac{3}{2})\gamma_{0}^{2}\omega_{d}<\alpha\wedge\frac{d}{2}$.
One can show, for all $x$, that $(\X_{0,\epsilon}(x))_{\epsilon >0}$
is a Cauchy sequence in $L^{2l}$ by computing
$E((\X_{0,\epsilon}^{j}(x)-\X_{0,\epsilon'}^{j}(x))^{2l})$ and letting
$\epsilon,\epsilon'$ go to $0$. Thus, $E((\X_{0}^{j}(x))^{2l})$ is given by
the left hand side of (\ref{eq:moment}).

To show the scaling (\ref{eq:scaling}), observe that we can prove the
following analogue to (\ref{eq:moment}) for any $q \leq 2l$:
 \begin{equation}\label{eq:momentincrement}
 E((\X_{0}^{j}(x+\lambda e)-\X_{0}^{j}(x))^{q})=
 \int_{(\R^{d})^{q}}\prod_{1 \leq i < j \leq
  q}\frac{e^{\gamma_{0}^{2}\phi(\frac{y_{i}-y_{j}}{R})}}{|\frac{y_{i}-y_{j}}{R}|_{*}^{\gamma_{0}^{2}\omega_{d}}}\prod_{1
  \leq i \leq q}f_{\lambda e}(y_{i})dy_{1} \ldots dy_{q},
\end{equation}   
where
\begin{equation}
f_{\lambda
  e}(y)=\frac{\varphi_{R}(y)y^{j}}{|y|^{d-\alpha+1}}-\frac{\varphi_{R}(y-\lambda e)(y^{j}-\lambda e^{j})}{|y-\lambda e|^{d-\alpha+1}}.
\end{equation}
By setting $y_{i}= \lambda u_{i}$ in the integral of
(\ref{eq:momentincrement}), we deduce easily (\ref{eq:scaling}).
\qed

\begin{rem}
It is not obvious why in the above proposition the coefficients
$C_{q}^{j}(e)$ are different from $0$ (cf. appendix).    
\end{rem}

\begin{rem}
Similarly as in the previous section, for $\gamma_{0}$ sufficiently
small, $\X_{0,\epsilon}$ converges in law to $\X_{0}$ in the space
of continuous fields. 
\end{rem}

\subsection{Nonsymmetry of the increments of $\X_{0}$}
By isotropy, we have:
\begin{equation*}
(\X_{0}(x+\lambda e)-\X_{0}(x)).e\underset{law}{=}\X_{0}^{j}(x+\lambda e_{j})-\X_{0}^{j}(x)
\end{equation*}

One can show that for $\lambda$ small the law is nonsymmetrical by
showing that the third moment is $\not= 0$, that is $C_{3}^{j}(e_{j}) \not=
0$ (see appendix).

\subsection{Towards a model of the turbulent velocity field}
In dimension 3, for the field $\X_{0}$, we have:
\begin{equation*}
\zeta_{q}=q\alpha-2 \pi q(q-1)\gamma_{0}^{2}. 
\end{equation*}
Thus, for $\alpha=1/3+4 \pi\gamma_{0}^{2}$, we have $\zeta_{3}=1$,
which means that for this choice the associated fields $\X_{0}$
satisfy at small scale the $4/5$ law with a non zero finite dissipation. 
Unfortunately, the fields in this family are not incompressible. The
incompressible case (at small scale) would correspond to the choice
$\alpha=1$, a limit case which is excluded by the constraint
$0<\alpha<1$ needed for the validity of the scaling of proposition
\ref{prop:scalingimp}. There is another severe obstacle for the choice
$\alpha=1$. Indeed, for the field $\X_{0}$ above, we have: 
\begin{equation*}
\zeta_{q}=(1/3+6\pi\gamma_{0}^{2})q-2 \pi\gamma_{0}^{2}q^{2}. 
\end{equation*}
One can easily identify the intermittency parameter
$4\pi\gamma_{0}^{2}$ using the experimental curve given in
\cite{cf:Anselmet} (cf. fig 8.8 p.132 in \cite{cf:Fri}). With this
data, we find $4\pi\gamma_{0}^{2}=0.023$. With this small
intermittency parameter, we would get $\alpha \sim 0.35$ which is not
close to the incompressible value $\alpha=1$. So, in spite of its
qualitative interest, this model cannot reach quantitative adequacy.

Another natural way to get incompressible fields is to use a
Biot-Savart like formula and take the limit as $\epsilon$ goes to $0$
of fields of the form:

\begin{equation*}
U^{\epsilon}(x)=\int_{\R^{3}}
 \varphi_{R}(x-y)\frac{x-y}{|x-y|_{\epsilon}^{d-\alpha+1}} \wedge
 d\Omega^{\epsilon},
 \end{equation*}
where $d\Omega^{\epsilon}$ is an isotropic random field. For example,
we can take: 
\begin{equation*}
d\Omega^{\epsilon}=e^{X^{\epsilon}(y)-C_{\epsilon}}dW(y),
\end{equation*}
where $dW(y)=(dW_{1}(y),dW_{2}(y),dW_{3}(y))$ denotes a three
dimensional white noise and $X^{\epsilon}$ is defined by the following formula:     
 \begin{equation*}
 X^{\epsilon}(y)=\gamma\int_{\R^{3}}K_{\epsilon}^{R}(y-\sigma).dW(\sigma),
 \end{equation*}
 with $K^{R}(x)=\frac{x}{|x|^{1+d/2}}1_{|x|\leq R}$.

As for the case of $\X_{0}$, we choose the constant $C_{\epsilon}$
such that $U^{\epsilon}$ converges to a non trivial field $U$ as
$\epsilon$ goes to zero. The vector field $U$ we obtain is
incompressible, homogeneous, isotropic and intermittent with
structural exponents $\zeta_{q}$ given by:
\begin{equation*}
\zeta_{q}=q\alpha-2\pi\gamma^{2}q(q-1).
\end{equation*}
Unfortunately, since the field $d\Omega^{\epsilon}(y)$ is isotropic
with respect to all unitary transformations (and not just the
rotations) we get for $U$ the symmetry:
\begin{equation*}
U(-x)-U(0)\underset{law}{=}U(x)-U(0)
\end{equation*}
so that the dissipation is equal to $0$. Thus the construction of an
homogeneous, isotropic, intermittent and incompressible vector field
with positive finite dissipation remains an open question.

\begin{comment}

In our approach, we perturb a Gaussian field to get multifractality and we further introduce some dependency to obtain also dissymmetry. This can make one think that in turbulence dissipation is linked to intermittency. This is a rather intricate issue.
On one hand, only dissymmetry seems to be needed to get energy dissipation (see the $4/5$ law).
On the other hand, it is well known experimentally that dissipation is not homogeneously distributed in the fluid but rather follows the lognormal distribution described by Kolmogorov and Obukhov by which it appears linked to intermittency.

In dimension $d=1$, our model displays some kind of (non causal) leverage effect. To get a realistic model for finance, with causal leverage effect, we have to make some specific change in the construction. This issue will be addressed in a forthcoming paper.

\end{comment}

\section{APPENDIX}

In this appendix, we prove that for $q$ even $C_{q}^{1}(e)$, given by
equation (\ref{eq:truccc}), is different from $0$ outside a
countable set and in the neighbourhood of $0$. We also show the same result for  $C_{3}^{j}(e_{j})$.

Consider first the case $q$ even; we set $q=2l$ with $l$ greater or equal to $1$ and we introduce the following function $F$: 
\begin{equation*}
F(\gamma)=\int_{(\R^{d})^{2l}}\prod_{1 \leq i < j \leq
  2l}\frac{1}{|u_{i}-u_{j}|^{\gamma}}\prod_{1
   \leq i \leq 2l}f(u_{i})du_{1} \ldots du_{2l},
\end{equation*}
where $f$ is the real function defined by: 
\begin{equation*}
f(u)=\frac{u^{j}+e^{j}/2}{|u+e/2|^{d-\alpha+1}}-\frac{u^{j}-e^{j}/2}{|u-e/2|^{d-\alpha+1}}.
\end{equation*}
The function $F$ is analytical in a neighbourhood of $0$; therefore,
in order to obtain the desired result, we have to prove that $F$ is
not identically equal to $0$. One can show that, for all $i < l$,
$F^{(i)}(0)=0$ and that:    
\begin{equation*}
F^{(l)}(0)=\frac{2l!}{2^{l}}(\int_{(\R^{d})^{2}}\ln(\frac{1}{|u_{1}-u_{2}|})f(u_{1})f(u_{2})du_{1}du_{2})^{l}.
\end{equation*}

The Fourier transform of $\ln(\frac{1}{|u|})$ is $a_{d}Pf(\frac{1}{|\xi|^d})+b_{d}\delta_{0}$ where $a_{d}>0$ and $b_{d}$ are two constants that depend only on the dimension and $Pf$ is Hadamard's finite part   (see p.258 in \cite{cf:Sch}). 
Since $\int_{\R^{d}}f(u)du=0$, we get:
\begin{equation*}
\int_{(\R^{d})^{2}}\ln(\frac{1}{|u_{1}-u_{2}|})f(u_{1})f(u_{2})du_{1}du_{2}=a_{d}\int_{\R^{d}}\frac{\hat{f}(\xi)^2}{|\xi|^d}d\xi,
\end{equation*}
thus $F^{(l)}(0)>0$.

Let us now consider $C_{3}^{j}(e_{j})$ and the corresponding function:
\begin{equation*}
F(\gamma)=\int_{(\R^{d})^3}\frac{1}{|u_{1}-u_{2}|^{\gamma}|u_{1}-u_{3}|^{\gamma}|u_{2}-u_{3}|^{\gamma}}f(u_{1})f(u_{2})f(u_{3})du_{1}du_{2}du_{3}.
\end{equation*}
We obviously have $F(0)=0$, $F'(0)=0$ and:
\begin{equation*}
\frac{1}{6}F''(0)=\mathcal{I}=\int_{(\R^{d})^3}\ln(|u_{1}-u_{2}|)\ln(|u_{1}-u_{3}|)f(u_{1})f(u_{2})f(u_{3})du_{1}du_{2}du_{3}
\end{equation*}
so that $\mathcal{I}=\int_{\R^{d}}f(x)\Theta(x)^2dx$ where:
\begin{equation*}
\Theta(x)=\int_{\R^{d}}\ln(|x-y|)f(y)dy.
\end{equation*}
Now we prove that there exists some real constant $c$ different from $0$ such that:  
\begin{equation}\label{eq:thetadef}
\Theta(x)=c(\frac{x^j+1/2}{|x+e_{j}/2|^{1-\alpha}}-\frac{x^j-1/2}{|x-e_{j}/2|^{1-\alpha}}).
\end{equation}

Indeed, we have (in what follows, $c$ denotes different real constants that are not equal to $0$):
\begin{equation*}
\hat{\Theta}(\xi)=c\frac{\hat{f}(\xi)}{|\xi|^d}
\end{equation*}
and $\hat{f}(\xi)=c\sin(\pi \xi^j)\frac{\xi^j}{|\xi|^{\alpha+1}}$ thus
\begin{equation*}
\hat{\Theta}(\xi)=c\sin(\pi \xi^j)\frac{\xi^j}{|\xi|^{d+\alpha+1}}.
\end{equation*}
The above expression (\ref{eq:thetadef}) now follows from:
\begin{equation*}
\widehat{(\frac{x^j}{|x|^{1-\alpha}})}(\xi)=c i \frac{\xi^j}{|\xi|^{d+\alpha+1}}
\end{equation*}
and $2i\sin(\pi \xi^j)=\widehat{\delta_{-e_{j}/2}}(\xi)-\widehat{\delta_{e_{j}/2}}(\xi)$.

Now let us denote $x=x^je_{j}+\tilde{x}$ and
\begin{equation*}
\phi(p,y,a)=\frac{y+1/2}{((y+1/2)^2+a)^p}-\frac{y-1/2}{((y-1/2)^2+a)^p}.
\end{equation*}
If we set:
\begin{equation*}
\varphi(x^j)=\phi(\frac{d-\alpha+1}{2},x^j,|\tilde{x}|^2)
\end{equation*}
and
\begin{equation*}
\psi(x^j)=\phi(\frac{1-\alpha}{2},x^j,|\tilde{x}|^2),
\end{equation*}
we get:
\begin{equation*}
\mathcal{I}=c^2\int_{\R^{d-1}}d\tilde{x}\int_{\R}\varphi(x^j)\psi(x^j)^2dx^j
\end{equation*}

Since $0<\alpha<1$, it is easy to check that $\psi(z)$ is a positive function of $z$, decreasing on $[0,\infty[$.
One can also check that there exists some $z^{*}>1/2$ such that $\varphi(z)$ is positive on $[0,z^{*}[$ and negative on $]z^{*},\infty[$.
Since $\varphi$ and $\psi$ are even and $\int_{0}^{\infty}\varphi(z)dz=0$, one can derive the following:
\begin{align*}
\int_{\R}\varphi(z)\psi(z)^2dz & = 2 \int_{0}^{\infty}\varphi(z)\psi(z)^2dz  \\
& =  2\int_{0}^{z^{*}}\varphi(z)\psi(z)^2dz+ 2\int_{z^{*}}^{\infty}\varphi(z)\psi(z)^2dz \\
& \geq  2\int_{0}^{z^{*}}\varphi(z)\psi(z)^2dz+ 2\psi(z^{*})^2\int_{z^{*}}^{\infty}\varphi(z)dz \\
& = 2\int_{0}^{z^{*}}\varphi(z)(\psi(z)^2-\psi(z^{*})^2)dz \\
& > 0. \\
\end{align*}
 It follows that $F''(0)>0$.

\qed



\end{document}